\documentclass[10pt]{article}
\usepackage{amssymb,amsmath,graphicx,acronym,amsfonts}

\newtheorem{corollary}{Corollary}[section]
\newtheorem{theorem}{Theorem}[section]
\newtheorem{lemma}{Lemma}[section]
\newtheorem{definition}{Definition}[section]
\newtheorem{proposition}{Proposition}[section]
\newtheorem{example}{Example}[section]
\newtheorem{assum}{Assumption}[section]
\newtheorem{algo}{Algorithm}[section]
\newtheorem{Remark}{Remark}[section]

\def\bc{\begin{corl}}
\def\bc{\end{corl}}
\def\ba{\begin{algo}}
\def\ea{\end{algo}}
\def\br{\begin{Remark}}
\def\er{\end{Remark}}
\def\bs{\begin{assum}}
\def\es{\end{assum}}
\def\bt{\begin{theorem}}
\def\et{\end{theorem}\vskip 3pt}
\def\bl{\begin{lemma}}
\def\el{\end{lemma}}
\def\ep{\end{proposition}}
\def\bp{\begin{proposition}}
\def\qed{\hfill{$\Box$}\vskip 5pt}
\def\be{\begin{example}}
\def\ee{\end{example}}
\def\bd{\begin{definition}}
\def\ed{\end{definition}}
\def\bc{\begin{corollary}}
\def\ec{\end{corollary}}
\def\proof{\noindent\it Proof. \hspace{1mm}\rm}
\input epsf
\begin{document}
\title{\bf Some Spectral Properties of Odd-Bipartite $Z$-Tensors and Their Absolute Tensors}
\author{Haibin Chen\thanks{Department of Applied Mathematics, The Hong Kong Polytechnic University, Hung Hom,
Kowloon, Hong Kong. Email: chenhaibin508@163.com. This author's work was supported by the Natural Science Foundation of China (11171180).}, \quad
Liqun Qi
\thanks{Department of Applied Mathematics, The Hong Kong Polytechnic University, Hung Hom,
Kowloon, Hong Kong. Email: maqilq@polyu.edu.hk. This author's work
was supported by the Hong Kong Research Grant Council (Grant No.
PolyU 502111, 501212, 501913 and 15302114).} }

\date{}
\maketitle
\begin{abstract}
Stimulated by odd-bipartite and even-bipartite hypergraphs, we define odd-bipartite (weakly odd-bipartie) and even-bipartite
(weakly even-bipartite) tensors. It is verified that all even order odd-bipartite tensors are irreducible tensors, while all
even-bipartite tensors are reducible no matter the parity of the order.
Based on properties of odd-bipartite tensors, we study the relationship between the largest H-eigenvalue of a $Z$-tensor with nonnegative diagonal elements, and the largest H-eigenvalue of absolute tensor of that $Z$-tensor. When the order is even and the $Z$-tensor is
weakly irreducible, we prove that the largest H-eigenvalue of the $Z$-tensor and the largest H-eigenvalue of the absolute tensor of that
$Z$-tensor are equal, if and only if the $Z$-tensor is weakly odd-bipartite.
Examples show the authenticity of the conclusions. Then,
we prove that a symmetric $Z$-tensor with nonnegative diagonal entries and the absolute tensor of the $Z$-tensor are diagonal similar, if and only if the $Z$-tensor has even order and it is weakly odd-bipartite.
After that, it is proved that, when an even order symmetric $Z$-tensor with nonnegative diagonal entries is weakly irreducible, the equality of the spectrum of the $Z$-tensor and the spectrum of absolute tensor of that $Z$-tensor, can be characterized by the equality of their spectral radii.
\medskip

\vskip 6pt


\end{abstract}

\newpage
\section{Introduction}

Since the early work of \cite{Lim05} and \cite{Qi05}, more and more researchers are interested in studying eigenvalue problems of tensors
in the past several years \cite{CS13, CPZ09, chen14,Chen14,Cooper12,Qi2013,Hu13,Ng09,Oeding13,Pearson13,QiL13,Qi13,Qi09,Qi08,Yang10,YY11}.
In \cite{Qi05}, two kinds of eigenvalues are defined for real symmetric tensors: eigenvalues
and E-eigenvalues. An eigenvalue (E-eigenvalue) with a real eigenvector (E-eigenvector)
is called an H-eigenvalue (Z-eigenvalue). When a symmetric tensor has even order, H-eigenvalues and
Z-eigenvalues always exist. An even order symmetric tensor is positive definite (semi-definite)
if and only if all of its H-eigenvalues or all of its Z-eigenvalues are positive (nonnegative).
Based upon this property, an H-eigenvalue method for the positive definiteness (positive semi-definiteness respectively) identification problem is developed.

The main difficulty in tensor problems is that they are generally nonlinear. Because of the difficulties in studying the properties of a general tensor, researchers focus on some structured tensors. $Z$-tensors are an important class of structured tensors and have been well studied \cite{Ding13,Mantica12,Zhang12}. They are closely related with spectral graph theory, the stationary distribution
of Markov chains and the convergence of iterative methods for linear systems.

Recently, in \cite{Hu13}, Hu et al. considered the largest Laplacian H-eigenvalue and the largest signless Laplacian
H-eigenvalue of a $k$-uniform connected hypergraph. When the order is even and the hypergraph is odd-bipartite, they proved that the largest Laplacian H-eigenvalue and the largest signless Laplacian H-eigenvalue are equal. For the odd order case,  it is
proved that the largest Laplacian H-eigenvalue is strictly less than the largest signless Laplacian H-eigenvalue \cite{Hu13}.
Later, Shao et al. \cite{Shao14} gave several spectral characterizations of the connected odd-bipartite hypergraphs.
They proved that the spectrum of the Laplacian tensor and the spectrun of the signless Laplacian tensor of an uniform hypergraph are equal if and only if the hypergraph is an even order connected odd-bipartite hypergraph.
Since the Laplacian tensor is a special case of
$Z$-tensors and the signless Laplacian tensor is a special case of the absolute tensors of $Z$-tensors, questions comes naturally: what is the relation between the largest H-eigenvalue of an general $Z$-tensor, and the largest H-eigenvalue of the $Z$-tensor's absolute tensor? What is the relation between spectrums of an general $Z$-tensor and its absolute tensor? These constitute main motivations of the paper.

In this article, some spectral properties of $Z$-tensors with nonnegative diagonal entries, and absolute tensors of $Z$-tensors
are studied. The rest of this paper is organized as follows.
In Section 2, some basic notions and preliminaries of tensors are presented.
In Section 3, stimulated by odd-bipartite and even-bipartite hypergraphs \cite{Qi2013}, odd-bipartite (weakly odd-bipartite) and even-bipartite (weakly even-bipartite) tensors are defined.
Odd-bipartite (even-bipartite) tensors are weakly odd-bipartite (weakly even-bipartite) tensors. Examples show that the converse, generally, may not hold.
An square odd-bipartite matrix is irreducible. For high order tensors,
we prove that an even order odd-bipartite tensor is irreducible,
while a tensor is reducible if it is even-bipartite no matter the parity of the order.

In Section 4, we study the relation between the largest H-eigenvalue of a $Z$-tensor with nonnegative diagonal entries, and the largest
H-eigenvalue of the $Z$-tensor's absolute tensor. For an even order $Z$-tensor with nonnegative diagonal entries, if it is weakly irreducible, we show that the largest H-eigenvalues of the $Z$-tensor and its absolute tensor are equal if and only if the $Z$-tensor is weakly odd-bipartite. For the odd order case, sufficient conditions for the equality of these largest H-eigenvalues are given. Examples show the authenticity of the conclusions.
In Section 5, we prove that, when an even order symmetric $Z$-tensor with nonnegative diagonal entries is weakly irreducible, its spectrum and the spectrum of its absolute tensor are equal if and only if the $Z$-tensor is odd-bipartite. Furthermore, it is shown that the equality of the spectrum of a symmetric $Z$-tensor with nonnegative diagonal entries, and the spectrum of the absolute tensor of that $Z$-tensor, can be characterized by the equality of their spectral radii. We conclude this paper with some final remarks in Section 6.

By the end of the introduction, we add some comments on notation that will be used in the sequel. Let $\mathbb{R}^n$
be the $n$ dimensional real Euclidean space and the set consisting of all
natural numbers is denoted by $\mathbb{N}$. Suppose $m, n\in \mathbb{N}$ are two natural numbers. Denote $[n]=\{1,2,\cdots,n\}$. Vectors are denoted by
italic lowercase letters i.e. $x,~ y,\cdots$, and tensors are written as calligraphic capitals such as
$\mathcal{A}, \mathcal{T}, \cdots.$ The $i$-th unit coordinate vector in $\mathbb{R}^n$ is denoted by $e_i$.
Let $|V|$ denote the number of elements when the symbol $|\cdot|$ be
used on a subset $V\subseteq \mathbb{N}$. If the symbol $|\cdot|$ is used on a tensor $\mathcal{A}=(a_{i_1 \cdots i_m})_{1\leq i_j\leq n}$, $j=1,\cdots,m$, we get another
tensor $|\mathcal{A}|=(|a_{i_1 \cdots i_m}|)_{1\leq i_j\leq n}$, $j=1,\cdots,m$.
If both $\mathcal{A}=(a_{i_1 \cdots i_m})_{1\leq i_j\leq n}$ and $\mathcal{B}=(b_{i_1 \cdots i_m})_{1\leq i_j\leq n}$, $j=1,\cdots,m$,
are real $m$th order $n$ dimensional tensors, then $\mathcal{A}\leq \mathcal{B}$ means $a_{i_1 \cdots i_m} \leq b_{i_1 \cdots i_m}$
for all $i_1,\cdots,i_m \in [n]$.

\setcounter{equation}{0}
\section{Preliminaries}
In this section, we will review some basic notions of tensors.
For more details, see \cite{Qi05} and the references therein.

A real $m$th order n-dimensional tensor $\mathcal{A}=(a_{i_1i_2\cdots i_m})$ is a multi-array
of real entries $a_{i_1i_2\cdots i_m}$, where $i_j \in [n]$ for $j\in [m]$.
If the entries $a_{i_1i_2\cdots i_m}$ are invariant under any permutation of
their indices, then tensor $\mathcal{A}$ is called a symmetric tensor.
The following definition on eigenvalue-eigenvector comes from \cite{Lim05, Qi05}.
\bd\label{def21} Let $\mathbb{C}$ be the complex field. A pair $(\lambda, x)\in \mathbb{C}\times \mathbb{C}^n\setminus \{0\}$ is called an
eigenvalue-eigenvector pair of $\mathcal{T}$, if they satisfy
\begin{equation}\label{e21}
\mathcal{T}x^{m-1}=\lambda x^{[m-1]},
\end{equation}
where $\mathcal{T}x^{m-1}$ and $x^{[m-1]}$ are all n dimensional column vectors such as
$$\mathcal{T}x^{m-1}=\left(\sum_{i_2,\cdots,i_m=1}^n t_{ii_2\cdots i_m}x_{i_2}\cdots x_{i_m} \right)_{1\leq i\leq n}$$ and $x^{[m-1]}=(x_i^{m-1})_{1\leq i\leq n}$.
\ed
For real tensor $\mathcal{T}$ and $x\in R^n$ in (\ref{e21}), $\lambda$ is a real number since
$\lambda=\frac{(\mathcal{T}x^{m-1})_j}{x^{m-1}_j}$
for some $j$ with $x_j\neq 0$. In this case,
$\lambda$ is called an H-eigenvalue of $\mathcal{T}$ and $x$ is its corresponding H-eigenvector \cite{Qi05}.

Next, we present a fundamental result which will be much used in the sequel.
\begin{proposition}\label{prop21}$^{\cite{Qi05}}$ Suppose that $\mathcal{T}=a(\mathcal{B}+b\mathcal{I})$, where $a$ and $b$ are two real numbers. Then $\mu$
is an eigenvalue (H-eigenvalue) of tensor $\mathcal{T}$ if and only if $\mu=a(\lambda+b)$, where $\lambda$ is an
eigenvalue (H-eigenvalue) of tensor $\mathcal{B}$. In this case, they have the same eigenvectors (H-eigenvectors).
\end{proposition}

The spectral radius of tensor $\mathcal{T}$ is denoted by
$$\rho(\mathcal{T})=\max \{|\lambda|: \lambda \text{ is an eigenvalue of}~ \mathcal{T}\}.$$
All eigenvalues of tensor $\mathcal{T}$ construct the spectrum denoted by $Spec(\mathcal{T})$.

\setcounter{equation}{0}
\section{Odd-bipartite and even-bipartite tensors}
In this section, we first define odd-bipartite tensors and even-bipartite tensors.
Then, some special characteristics of this kinds of tensors are shown.
\bd\label{def41} Assume $\mathcal{A}=(a_{i_1\cdots i_m})$ is an tensor with order m and dimension n. If there is a nonempty proper
index subset $V\subset [n]$ such that
$$a_{i_1\cdots i_m}\neq 0, \text{ when $|V\cap \{i_1,\cdots,i_m\}|$ is odd}$$
and $a_{i_1\cdots i_m}=0$ for the others, then $\mathcal{A}$ is called an odd-bipartite tensor corresponding to set $V$ or $\mathcal{A}$ is odd-bipartite
for simple.
\ed
Here, we should note that indices of an edge $\{i_1,\cdots,i_m\}$ in hypergraph \cite{Qi2013} are different from each other,
which is a notable distinction to general tensors. So, in this paper, we define that $|V\cap \{i_1,\cdots,i_m\}|$ is the number of
indices $V\cap \{i_1,\cdots,i_m\}$, and duplicate
indices should be calculated. For example, suppose $V=\{1,2,3\}$ and $\mathcal{A}$ is a 4th order 6 dimensional tensor, then
$$|V\cap \{1,1,3,3\}|=4,~|V\cap \{1,2,3,5\}|=3, ~|V\cap \{4,6,4,5\}|=0.$$
\bd\label{def42} Assume $\mathcal{A}=(a_{i_1\cdots i_m})$ is a tensor with order $m$ and dimension $n$. $\mathcal{A}$
is called weakly odd-bipartite if there is a nonempty
proper index subset $V\subset [n]$ such that
$$a_{i_1\cdots i_m}= 0, \text{ when $|V\cap \{i_1,\cdots,i_m\}|$ is even}.$$
\ed

From Definitions \ref{def41} and \ref{def42}, even-bipartite and weakly even-bipartite tensors can be defined similarly. Furthermore, we can easily prove that, if $\mathcal{A}$ is odd-bipartite (even-bipartite, respectively), then $\mathcal{A}$ is
weakly odd-bipartite (weakly even-bipartite respectively), but not vice versa. For example, suppose $\mathcal{A}$ is a $3$th order 2 dimensional tensor with
entries such that
$$a_{222}=1~\text{and}~a_{i_1i_2i_3}=0$$
for the others. It is easy to check that $\mathcal{A}$ is weakly odd-bipartite corresponding to
the index set $V=\{2\}$ but not odd-bipartite corresponding to $\{1\}$ or $\{2\}$.

When $m$ is odd, for all $i_1,i_2,\cdots,i_m\in [n]$ and a nonempty proper index subset $V\subset [n]$, it holds that
$|\{i_1,i_2,\cdots,i_m\}\cap V|$ is odd if and only if $|\{i_1,i_2,\cdots,i_m\}\cap \bar{V}|$ is even, where $\bar{V}=[n]\backslash V$.
So, by Definitions \ref{def41} and \ref{def42}, it's easy to check that the following result holds.
\bl\label{lema41} Let $\mathcal{A}$ be a tensor with order $m$ and dimension $n$. Assume $m$ is odd. Then, $\mathcal{A}$
is odd-bipartite (or weakly odd-bipartite respectively) corresponding to nonempty proper index subset $V\subset [n]$ if and only if
$\mathcal{A}$ is even-bipartite (or weakly even-bipartite respectively) corresponding to the nonempty proper index subset $\bar{V}=[n]\backslash V$.
\el

Irreducible tensors are a class of important and useful tensors, which have been repeatedly used in Perron Frobenius Theorem for nonnegative tensors \cite{Chang08, Yang10, YY11}. Next, we will study
the relation between irreducible tensors and odd-bipartite tensors. To do this, we first list the corresponding definition below.
\bd\label{def43}$^{\cite{Chang08}}$ For a tensor $\mathcal{T}$ with order m and dimension n. We call $\mathcal{T}$ is reducible if there is a nonempty
proper index subset $V\subset [n]$ such that
$$t_{i_1i_2\cdots i_m}=0,~\forall ~i_1\in V,~\forall~i_2,i_3,\cdots,i_m \notin V.$$
Otherwise we call $\mathcal{T}$ is irreducible.
\ed

\bt\label{them41} Let $m$ be even. Assume tensor $\mathcal{A}=(a_{i_1\cdots i_m})$ with order m and dimension n is odd-bipartite. Then
$\mathcal{A}$ is irreducible.
\et

\proof Since $\mathcal{A}$ is odd-bipartite, there exists a nonempty proper index subset $V\subset [n]$ satisfying

\begin{equation}\label{e41}
a_{i_1\cdots i_m}\neq 0, \text{ when the number $|V\cap \{i_1,\cdots,i_m\}|$ is odd}.
\end{equation}

By contradiction, suppose $\mathcal{A}=(a_{i_1\cdots i_m})$ is reducible, then there is a nonempty
proper index subset $V_1\subset [n]$ such that
\begin{equation}\label{e42}
a_{i_1\cdots i_m}=0,~\forall~i_1\in V_1,~\forall~i_2,\cdots,i_m \notin V_1.
\end{equation}

(i) If $V_1\subseteq V$, let $i_1\in V_1, i_2,\cdots,i_m\notin V$. Here, several indices in $i_2,\cdots, i_m$ may equal to each other when the number
of elements in $[n]\backslash V$ is strictly less than $m-1$.
Then, by (\ref{e42}) we have
$$a_{i_1\cdots i_m}=0,$$
which contradicts with (\ref{e41}) since $|V\cap \{i_1,\cdots,i_m\}|=1$ is odd.

(ii) If $V\subseteq V_1$, let $i_1\in V, i_2,\cdots,i_m\notin V_1$. Then, by (\ref{e42}) one has
$$a_{i_1\cdots i_m}=0,$$
which is a contradiction with (\ref{e41}).

(iii) If $V\cap V_1\neq \emptyset$ and neither $V\subseteq V_1$ nor $V_1\subseteq V$, let $i_1\in V_1 \setminus V, ~i_2,~\cdots,i_m\in V\setminus V_1$. Then it follows that
$$a_{i_1\cdots i_m}=0,$$
which also contradicts (\ref{e41}), since $|V\cap \{i_1,\cdots,i_m\}|=m-1$ is a odd number.

(iv) If $V\cap V_1=\emptyset$, let $i_1\in V_1, i_2,\cdots,i_m \in V$.
By Definition \ref{def43}, we have
$$a_{i_1\cdots i_m}=0.$$
Since $|V\cap \{i_1,\cdots,i_m\}|=m-1$ be odd, by (\ref{e41}), one has
$$a_{i_1\cdots i_m}\neq 0,$$
which is a contradiction.
All in all, we know that $\mathcal{A}$ can not be reducible and the desired results follows. \qed

If a tensor $\mathcal{A}$ is even-bipartite, no matter the order of $\mathcal{A}$ is odd or even, we have the following result.
\bt\label{them42} Assume tensor $\mathcal{A}=(a_{i_1\cdots i_m})$ with order m and dimension n is even-bipartite corresponding to
a nonempty proper index subset $V\subseteq [n]$. Then
$\mathcal{A}$ is reducible corresponding to $V$.
\et
\proof By  definitions of reducible tensors and even-bipartite tensors, the conclusion obviously holds. \qed

Suppose an even order $Z$-tensor and its absolute tensor are defined such that,
\begin{equation}\label{e43} \mathcal{A}=\mathcal{D}-\mathcal{C},~~|\mathcal{A}|=\mathcal{D}+\mathcal{C},\end{equation}
where $\mathcal{D}$ is an nonnegative diagonal tensor and $\mathcal{C}$ is an nonnegative tensor.
From Theorem \ref{them42}, if  $\mathcal{C}$ is odd-bipartite,
then tensors $\mathcal{A}$ and $|\mathcal{A}|$ are irreducible. Combining this with Theorem 3.1 of \cite{Fried13} we have the following result.
\bc\label{corl41} Let $m$ be even. Suppose tensor $\mathcal{A}=\mathcal{D}-\mathcal{C}$ with order m and dimension n is defined as in (\ref{e43}).
Then, $\mathcal{A}$ and its absolute tensor $|\mathcal{A}|$ are all weakly irreducible if nonnegative tensor $\mathcal{C}$ is odd-bipartite.
\ec

By the Perron-Frobenius theorem on nonnegative tensors in \cite{Chang08} and by Theorem 4.1 of \cite{Fried13}, the following result follows.

\bc\label{corl42} Let $m$ be even. Assume tensor $\mathcal{A}$ is defined as in Corollary \ref{corl41}. If $\mathcal{C}$ is odd-bipartite,
the largest H-eigenvalue of $|\mathcal{A}|$ is $\rho(|\mathcal{A}|)$.
Furthermore, there exists a positive $n$ dimensional eigenvector $x\in \mathbb{R}^n$ such that
$$|\mathcal{A}|x^{m-1}=\rho(|\mathcal{A}|)x^{[m-1]}.$$
\ec

\setcounter{equation}{0}
\section{The relation between the largest H-eigenvalues of a $Z$-tensor and its absolute tensor}
In this section, suppose an order $m$ dimension $n$ $Z$-tensor $\mathcal{A}$ with nonnegative diagonal elements has format
\begin{equation}\label{e51}
\mathcal{A}=\mathcal{D}-\mathcal{C},
\end{equation}
where $\mathcal{D}$ is an nonnegative diagonal tensor and $\mathcal{C}$ is an nonnegative tensor with zero diagonal elements.
So the absolute format of $\mathcal{A}$ is $|\mathcal{A}|=\mathcal{D}+\mathcal{C}.$
In the following analysis, entries of $\mathcal{A},~\mathcal{C}$ and $\mathcal{D}$ are always defined as below
$$\mathcal{A}=(a_{i_1\cdots i_m}),~\mathcal{C}=(c_{i_1\cdots i_m}), ~\mathcal{D}=(d_{i_1\cdots i_m}),~i_1,i_2,\cdots,i_m\in [n].$$
For the sake of simple, let $d_{ii\cdots i}=d_i$, $i\in [n]$.

During this part, we mainly study the relationship between the
largest H-eigenvalue of a $Z$-tensor $\mathcal{A}$ in (\ref{e51}),
and the largest H-eigenvalue of the absolute tensor of $\mathcal{A}$. Sufficient and necessary conditions or sufficient conditions to
guarantee the equality of these largest H-eigenvalues are shown. It should be noted that
all even order nonnegative tensors always have H-eigenvalues \cite{Yang10}. To proceed, we make an assumption in advance,
all tensors considered in this part always have H-eigenvalues.

The largest H-eigenvalues of $\mathcal{A}$ and $|\mathcal{A}|$
are denoted by $\lambda(\mathcal{A})$ and $\lambda(|A|)$ respectively. From Corollary \ref{corl42}, we know that $\lambda(|A|)=\rho(|\mathcal{A}|)$.

\bt\label{thm51} Let $m$ be even. Suppose $\mathcal{A}=\mathcal{D}-\mathcal{C}$ is defined as (\ref{e51}). Then,
$$\lambda(\mathcal{A})=\lambda(|A|)$$
if $\mathcal{C}$ is odd-bipartite.
\et

\proof
By Lemma 13 of \cite{Qi13}, we have
$$\lambda(\mathcal{A})\leq \rho(\mathcal{A})\leq \rho(|\mathcal{A}|)=\lambda(|\mathcal{A}|).$$
Thus, in order to prove the conclusion, we only need to prove
$$\lambda(|\mathcal{A}|)\leq\lambda(\mathcal{A}).$$
Since $\mathcal{C}$ is odd-bipartite, there exists a nonempty proper index subset $V\subset [n]$ satisfying
$$c_{i_1\cdots i_m}\neq 0,~\text{if $|V\cap \{i_1,\cdots,i_m\}|$ is odd,}$$
and $c_{i_1\cdots i_m}=0$ for the others. So, for all entries of $\mathcal{A}$, it follows that
$$a_{i_1\cdots i_m}\neq 0,~\text{if $|V\cap \{i_1,\cdots,i_m\}|$ is odd,}$$
and $a_{i_1\cdots i_m}=0$ for the others except the diagonal entries $a_{ii\cdots i},~i\in [n]$.
By Theorem \ref{them42}, we know that $\mathcal{C}$, $\mathcal{A}$ and $|\mathcal{A}|$ are all irreducible tensors.
From Theorem 4.1 of \cite{Fried13} and Definition
\ref{def21}, there is a vector $x\in \mathbb{R}^n$, $x>{\bf 0}$ satisfying
$$|\mathcal{A}|x^{m-1}=\lambda(|\mathcal{A}|)x^{[m-1]}.$$
Suppose $y\in \mathbb{R}^n$ be defined such that $y_i =x_i$ whenever $i\in V$
and $y_i=-x_i$ for the others. When $i\in V$, we have

\begin{equation}\label{e52}
\begin{array}{rl}(\mathcal{A}y^{m-1})_i=&\left[(\mathcal{D}-\mathcal{C})y^{m-1}\right]_i\\
=&d_iy_i^{m-1}-\sum_{i_2,\cdots,i_m \in [n]}c_{ii_2\cdots i_m}y_{i_2} \cdots y_{i_m}\\
=&d_iy_i^{m-1}-\sum_{i_2,\cdots,i_m \in [n]~\text{ $|V\cap \{i,i_2,\cdots, i_m\}|$ is odd}}c_{ii_2\cdots i_m}y_{i_2} \cdots y_{i_m}\\
=&d_ix_i^{m-1}+\sum_{i_2,\cdots,i_m \in [n]~\text{ $|V\cap \{i,i_2,\cdots, i_m\}|$ is odd}}c_{ii_2\cdots i_m}x_{i_2} \cdots x_{i_m}\\
=&\left[(\mathcal{D}+\mathcal{C})x^{m-1} \right]_i\\
=&\lambda(|\mathcal{A}|)x_i^{m-1}\\
=&\lambda(|\mathcal{A}|)y_i^{m-1}.
\end{array}
\end{equation}
Here the fourth equality follows the fact that $m$ is even and exactly odd number indices take negative values for
each $\{i_2,\cdots,i_m\}\subseteq [n]$. When $i\notin V$, we have
\begin{equation}\label{e53}
\begin{array}{rl}(\mathcal{A}y^{m-1})_i=&\left[(\mathcal{D}-\mathcal{C})y^{m-1}\right]_i\\
=&d_iy_i^{m-1}-\sum_{i_2,\cdots,i_m \in [n]}c_{ii_2\cdots i_m}y_{i_2} \cdots y_{i_m}\\
=&d_iy_i^{m-1}-\sum_{i_2,\cdots,i_m \in [n]~\text{ $|V\cap \{i,i_2,\cdots, i_m\}|$ is odd}}c_{ii_2\cdots i_m}y_{i_2} \cdots y_{i_m}\\
=&-d_ix_i^{m-1}-\sum_{i_2,\cdots,i_m \in [n]~\text{ $|V\cap \{i,i_2,\cdots, i_m\}|$ is odd}}c_{ii_2\cdots i_m}x_{i_2} \cdots x_{i_m}\\
=&-\left[(\mathcal{D}+\mathcal{C})x^{m-1}\right]_i\\
=&-\lambda(|\mathcal{A}|)x_i^{m-1}\\
=&\lambda(|\mathcal{A}|)y_i^{m-1}.
\end{array}
\end{equation}
Here the fourth equality follows the fact that $m$ is even and exactly even number indices take negative values for
each $\{i_2,\cdots,i_m\}\subseteq [n]$. The last equality of (\ref{e53})
follows from the definition of $y_i=-x_i$ when $i\notin V$. Thus, by (\ref{e52}), (\ref{e53}) and
Definition \ref{def21}, $\lambda(|\mathcal{A}|)$ is a H-eigenvalue of $\mathcal{A}$ with H-eigenvector $y$. So, we have
$$\lambda(|\mathcal{A}|)\leq\lambda(\mathcal{A}),$$
and the desired result follows. \qed

Here, in the proof of Theorem \ref{thm51}, odd-bipartite property of $\mathcal{C}$ guarantees that $|\mathcal{A}|$ has a positive H-eigenvector.
Actually, if the H-eigenvector is nonnegative, one can obtain the same result. Before proving this, we first cite an useful conclusion.

\bl\label{lema52}$^{\cite{Yang10}}$ If $\mathcal{A}$ is a nonnegative tensor with order $m$ and dimension $n$, then $\rho(\mathcal{A})$ is an eigenvalue of $\mathcal{A}$ with
a nonnegative eigenvector $y\neq0$.
\el

\bt\label{thm52} Let $m$ be even. Suppose $\mathcal{A}$ is defined as in Theorem \ref{thm51}. If $\mathcal{C}$ is weakly odd-bipartite,
 then it holds that
$$\lambda(\mathcal{A})=\lambda(|A|).$$
\et
\proof Since tensor $\mathcal{C}$ is weakly odd-bipartite, so there is a nonempty proper index subset $V\subseteq [n]$ such that
$$c_{i_1\cdots i_m}=0,~\text{when}~|\{i_1,\cdots, i_m\}\cap V|~ \text{is even},$$
and
$|\{i_1,\cdots, i_m\}\cap V|$ must be an odd number for nonzero entries $c_{i_1\cdots i_m}\neq0$, $i_1,\cdots, i_m \in [n]$.

On the other hand, by Lemma \ref{lema52}, there is a nonnegative H-eigenvector $x\geq0$ of $|\mathcal{A}|$ corresponding
to $\lambda(|\mathcal{A}|)$ .
Suppose vector $y\in \mathbb{R}^n$ be defined such that $y_i =x_i$ whenever $i\in V$ and $y_i=-x_i$ for the others.
Then, the remaining process is similar with the proof of Theorem \ref{thm51}. \qed

Now, we will give an example to show that the conditions in Theorem \ref{thm52} is not necessary. For example, suppose $4$th order 2 dimensional tensor $\mathcal{A}$ with entries such that
$$a_{1111}=a_{2222}=1,~a_{1122}=-1,$$
and $a_{i_1i_2i_3i_4}=0$ for the others.
After calculating the largest H-eigenvalues of $\mathcal{A}$ and $|\mathcal{A}|$, we obtain
$$\lambda(\mathcal{A})=\lambda(|\mathcal{A}|)=1.$$
But, the nonnegative tensor $\mathcal{C}$ is not weakly odd-bipartite corresponding to any nonempty proper index subset of $\{1,2\}$.
In the following, sufficient and necessary conditions for the equality of the two largest
H-eigenvalues are presented, and it is proved that the necessity of the Theorem \ref{thm52} holds when the nonnegative tensor $\mathcal{C}$ is weakly irreducible. Before doing this, we cite a definition.

\bd\label{def51}$^{\cite{Qi13}}$ Assume that $\mathcal{T}$ is a tensor with order m and dimension n. Construct a graph
$\hat{G}=(\hat{V}, \hat{E})$, where $\hat{V}=\cup_{j=1}^d V_j$ and $V_j$ are subsets of $\{1,2,\cdots,n\}$
for $j=1,\cdots,d$. Suppose that $i_j\in V_j, i_l\in V_l, j\neq l$. $(i_j, i_l)\in \hat{E}$ if and only if
$t_{i_1i_2\cdots i_m}\neq0$ for some $m-2$ indices $\{i_1,\cdots,i_m \} \backslash \{i_j, i_l\}$. Then, tensor
$\mathcal{T}$ is called weakly irreducible if $\hat{G}$ is connected.
\ed
As observed in \cite{Fried13}, an irreducible tensor must be always weakly irreducible.

\bt\label{thm53} Let $\mathcal{A}$ be defined as in Theorem \ref{thm52}. Assume $\mathcal{C}$ is weakly irreducible. Then,
$$\lambda(\mathcal{A})=\lambda(|A|),$$
if and only if $\mathcal{C}$ is weakly odd-bipartite.
\et
\proof The sufficient condition has been proved in Theorem \ref{thm52}, and we only need to prove the necessary part.

Suppose $x\in \mathbb{R}^n$ is an H-eigenvector of $\mathcal{A}$ corresponding to $\lambda(\mathcal{A})$ such that
$\sum_{i=1}^{n}x_i^m=1$. Assume $y\in \mathbb{R}^n$ be defined by $y_i=|x_i|$, for $i\in [n]$. Since $m$ is even,
one has $\sum_{i=1}^{n}y_i^m=1$.
By Lemma 3.1 of \cite{Li13}, we have
\begin{equation}\label{e54}
\begin{array}{rl} \lambda(\mathcal{A})=&\mathcal{A}x^m=(\mathcal{D}-\mathcal{C})x^m\\
=&\sum_{i=1}^{n}d_ix_i^m-\sum_{i_1,\cdots,i_m \in [n]}c_{i_1i_2\cdots i_m}x_{i_1} \cdots x_{i_m}\\
\leq & \sum_{i=1}^{n}d_iy_i^m+\sum_{i_1,\cdots,i_m \in [n]}c_{i_1i_2\cdots i_m}y_{i_1} \cdots y_{i_m}\\
=&(\mathcal{D}+\mathcal{C})y^m\leq \lambda(|\mathcal{A}|).
\end{array}
\end{equation}
Hence, by the fact that $\lambda(\mathcal{A})=\lambda(|\mathcal{A}|)$, all inequalities in equation (\ref{e54}) should be equalities, which implies that
$y$ is a H-eigenvector of $|\mathcal{A}|$ corresponding to $\lambda(|\mathcal{A}|)$. Since $\mathcal{C}$ is weakly irreducible,
$|\mathcal{A}|$ is also weakly irreducible. According to Theorem 4.1 of \cite{Fried13},
it holds that $y>0$ i.e., all elements in $y$ are positive. Let $V=\{i\in [n]|~x_i>0 \}$ and $\bar{V}=\{i\in [n]|~x_i<0 \}$.
Then $V\cup \bar{V}=[n]$. By (\ref{e54}), we obtain
$$\sum_{i_1,\cdots,i_m \in [n]}c_{i_1i_2\cdots i_m}(|x_{i_1}| \cdots |x_{i_m}|+x_{i_1} \cdots x_{i_m})=0,$$
which implies that
$$c_{i_1i_2\cdots i_m}(|x_{i_1}| \cdots |x_{i_m}|+x_{i_1} \cdots x_{i_m})=0,
$$
for all $i_1,i_2,\cdots,i_m \in [n]$ since $\mathcal{C}$ is nonnegative. When $|\{i_1,i_2,\cdots,i_m\}\cap V|$ is even, we have
$$|x_{i_1}| \cdots |x_{i_m}|+x_{i_1} \cdots x_{i_m}>0,$$
which implies $c_{i_1i_2\cdots i_m}=0$.
When $|\{i_1,i_2,\cdots,i_m\}\cap V|$ is odd, we have
$$|x_{i_1}| \cdots |x_{i_m}|+x_{i_1} \cdots x_{i_m}=0.$$
In this case, the value $c_{i_1i_2\cdots i_m}$ may be zero or may not be zero . Thus,
from Definition \ref{def42}, it follows that $\mathcal{C}$
is weakly odd-bipartite corresponding to set $V$ and the desired conclusion holds. \qed

Next, we study the relationship between a $Z$-tensor and its absolute tensor in the odd order case.
In \cite{Hu13}, Hu et al. proved that the largest H-eigenvalue of an odd order Laplacian tensor is always strictly less than
the largest H-eigenvalue of an signless Laplacian tensor corresponded to the Laplacian tensor. By definitions of Laplacian tensor and signless Laplacian tensor in connected hypergraphs, we know that their diagonal entries
are positive, and subscripts of each nonzero element are mutually distinct.
However, general $Z$-tensors (\ref{e51}) may not possess those advantages.
Hence, for a general odd order $Z$-tensor (\ref{e51}), the largest H-eigenvalue of $\mathcal{A}$ may not be strictly less
than the largest H-eigenvalue of $|\mathcal{A}|$ when the order is odd.

The following example shows that the largest H-eigenvalues of a $Z$-tensor (\ref{e51}) and its absolute tensor are equal.
\begin{example}\label{exam51} Let $\mathcal{A}$ be a 5th order 3 dimensional tensor. Its entries are given by
$$a_{11111}=a_{22222}=a_{33333}=1,~a_{11122}=a_{22233}=-1$$
and $a_{i_1i_2i_3i_4i_5}=0$ for the others. Then the H-eigenvalue problems for $\mathcal{A}$ and $|\mathcal{A}|$ are
$$
\begin{cases}
x^4_1-x^2_1x^2_2=\lambda x^4_1\\
x^4_2-x_2^2x_3^2=\lambda x_2^4\\
x_3^4=\lambda x_3^4\\
\end{cases}
$$
and
$$
\begin{cases}
x^4_1+x^2_1x^2_2=\lambda x^4_1\\
x^4_2+x_2^2x_3^2=\lambda x_2^4\\
x_3^4=\lambda x_3^4\\
\end{cases}.
$$
After calculating these equation sets, we know that $\lambda(\mathcal{A})=\lambda(|\mathcal{A}|)=1$.
\end{example} \qed

\bt\label{thm54} Let $A$ be defined as (\ref{e51}). Assume $m$ is odd. Suppose $\mathcal{C}$ is
weakly odd-bipartite corresponding to a nonempty proper index subset $V\subseteq [n]$. If for all $i\in V$, it satisfies
$$c_{ii_2i_3\cdots i_m}=0,~\forall ~i_2,i_3,\cdots, i_m\in [n],$$
then $\lambda(\mathcal{A})=\lambda(|\mathcal{A}|).$
\et
\proof By the analysis in Theorems \ref{thm51}-\ref{thm53}, from Lemma 13 of \cite{Qi13} and Corollary \ref{corl42},
it follows that
$$\lambda(\mathcal{A})\leq \rho(\mathcal{A})\leq \rho(|\mathcal{A}|)= \lambda(|\mathcal{A}|).$$
Thus, we only need to prove $$\lambda(|\mathcal{A}|)\leq \lambda(\mathcal{A}).$$
Let $x\in \mathbb{R}^n$ be a nonnegative H-eigenvector of $|\mathcal{A}|$ corresponding to $\lambda(|\mathcal{A}|)$.
So, for all $i\in [n]$, we have
\begin{equation}\label{e55}
(|\mathcal{A}|x^{m-1})_i=[(\mathcal{D}+\mathcal{C})x^{m-1}]_i=\lambda(|\mathcal{A}|)x_i^{m-1}.
\end{equation}
Suppose $y\in \mathbb{R}^n$ be defined as $y_i=-x_i$, $i\in V$ and $y_i=x_i$, $i\notin V$.
By conditions, $\mathcal{C}$ is weakly odd-bipartite corresponding to subset $V$, which means
$$c_{i_1i_2i_3\cdots i_m}=0,~i_1,i_2,\cdots, i_m\in [n]$$
when $|\{i_1,i_2,i_3,\cdots, i_m\}\cap V|$ is even.
Then, for all $i\in [n]$, one has
\begin{equation}\label{e56}
\begin{array}{rl} (\mathcal{A}y^{m-1})_i=&[(\mathcal{D}-\mathcal{C})y^{m-1}]_i\\
=&d_iy_i^{m-1}-\sum_{i_2,\cdots,i_m \in [n]~\text{ $|V\cap \{i,i_2,\cdots, i_m\}|$ is odd}}c_{ii_2\cdots i_m}y_{i_2} \cdots y_{i_m}\\
=&d_ix_i^{m-1}-\sum_{i_2,\cdots,i_m \in [n]~\text{ $|V\cap \{i,i_2,\cdots, i_m\}|$ is odd}}c_{ii_2\cdots i_m}y_{i_2} \cdots y_{i_m},\\
\end{array}
\end{equation}
where the third equality follows $m-1$ is even and $y_i^{m-1}=x_i^{m-1}$.
When $i\in V$, by the fact that $c_{ii_2i_3\cdots i_m}=0,~i_2,i_3,\cdots, i_m\in [n]$, and by (\ref{e55}), (\ref{e56}), we have
\begin{equation}\label{e57}
\begin{array}{rl} (\mathcal{A}y^{m-1})_i=&[(\mathcal{D}-\mathcal{C})y^{m-1}]_i\\
=&d_iy_i^{m-1}-\sum_{i_2,\cdots,i_m \in [n]}c_{ii_2\cdots i_m}y_{i_2} \cdots y_{i_m}\\
=&d_ix_i^{m-1}=\lambda(|\mathcal{A}|)x_i^{m-1}\\
=&\lambda(|\mathcal{A}|)y_i^{m-1}.
\end{array}
\end{equation}
Similarly, when $i\notin V$, it holds that
\begin{equation}\label{e58}
\begin{array}{rl} (\mathcal{A}y^{m-1})_i=&[(\mathcal{D}-\mathcal{C})y^{m-1}]_i\\
=&d_iy_i^{m-1}-\sum_{i_2,\cdots,i_m \in [n]~\text{ $|V\cap \{i,i_2,\cdots, i_m\}|$ is odd}}c_{ii_2\cdots i_m}y_{i_2} \cdots y_{i_m}\\
=&d_ix_i^{m-1}+\sum_{i_2,\cdots,i_m \in [n]~\text{ $|V\cap \{i,i_2,\cdots, i_m\}|$ is odd}}c_{ii_2\cdots i_m}x_{i_2} \cdots x_{i_m}\\
=&d_ix_i^{m-1}+(\mathcal{C}x^{m-1})_i\\
=&[(\mathcal{D}+\mathcal{C})x^{m-1}]_i\\
=&(|\mathcal{A}|x^{m-1})_i=\lambda(|\mathcal{A}|)x_i^{m-1}\\
=&\lambda(|\mathcal{A}|)y_i^{m-1},
\end{array}
\end{equation}
where the third equality follows the fact that $m$ is odd and exactly odd indices take negative values.
By (\ref{e57}) and (\ref{e58}), we know that $\lambda(|\mathcal{A}|)$ is a H-eigenvalue of $\mathcal{A}$.
Hence, we have $\lambda(|\mathcal{A}|)\leq \lambda(\mathcal{A})$ and the desired result follows.
\qed
Now, we present a example to verify the authenticity of Theorem \ref{thm54}.
\begin{example}\label{exam52} Set a 5th order 3 dimensional tensor $\mathcal{A}$ such that
$$a_{11111}=1,~a_{22222}=1,~a_{33333}=3,~a_{11333}=-1,~a_{22333}=-2$$
and $a_{i_1i_2i_3i_4}=0$ for the others. Let $V=\{3\}$.
Then $\mathcal{C}$ is weakly odd-bipartite corresponding to the set $V$ and $c_{3i_2i_3i_4i_5}=0,~ \forall~ i_2,i_3,i_4,i_5\in [3]$.

The H-eigenvalue problems for $\mathcal{A}$ and $|\mathcal{A}|$
are to solve
$$
\begin{cases}
x^4_1-x_1x_3^3=\lambda x^4_1\\
x^4_2-2x_2x_3^3=\lambda x_2^4\\
3x_3^4=\lambda x_3^4
\end{cases}
$$
and
$$
\begin{cases}
x^4_1+x_1x_3^3=\lambda x^4_1\\
x^4_2+2x_2x_3^3=\lambda x_2^4\\
3x_3^4=\lambda x_3^4
\end{cases}.
$$
After calculating the largest H-eigenvalues of $\mathcal{A}$ and $|\mathcal{A}|$, we obtain
$$\lambda(\mathcal{A})=\lambda(|\mathcal{A}|)=3.$$ \qed
\end{example}
The next example shows that the conditions in Theorem \ref{thm54} are not necessary.
\begin{example}\label{exam53} Let $\mathcal{A}$ be a 5th order 3 dimensional tensor. Its entries are given by
$$a_{11111}=1,~a_{22222}=2,~a_{33333}=4,~a_{11122}=a_{11333}=-1,~a_{22233}=-2$$
and $a_{i_1i_2i_3i_4i_5}=0$ for the others. Then the H-eigenvalue problems for $\mathcal{A}$ and $|\mathcal{A}|$ are
$$
\begin{cases}
x^4_1-x^2_1x^2_2-x_1x_3^3=\lambda x^4_1\\
2x^4_2-2x_2^2x_3^2=\lambda x_2^4\\
4x_3^4=\lambda x_3^4\\
\end{cases}
$$
and
$$
\begin{cases}
x^4_1+x^2_1x^2_2+x_1x_3^3=\lambda x^4_1\\
2x^4_2+2x_2^2x_3^2=\lambda x_2^4\\
4x_3^4=\lambda x_3^4\\
\end{cases}
$$
After calculating these equation sets, we know that $\lambda(\mathcal{A})=\lambda(|\mathcal{A}|)=4$, but the
nonnegative tensor $\mathcal{C}$ is not weakly odd-bipartite corresponding to any nonempty proper index subset of $\{1,2,3\}$.
\end{example} \qed

By Lemma \ref{lema41} and Theorem \ref{thm54}, we have the following conclusion.
\bc\label{corl51}Let $A$ be defined as in (\ref{e51}). Assume $m$ is odd. Suppose $\mathcal{C}$ is
weakly even-bipartite corresponding to a nonempty proper index subset $V\subseteq [n]$. If for all $i\notin V$, it satisfies
$$c_{ii_2i_3\cdots i_m}=0,~\forall ~i_2,i_3,\cdots, i_m\in [n],$$
then $\lambda(\mathcal{A})=\lambda(|\mathcal{A}|).$
\ec

\setcounter{equation}{0}
\section{The relation between spectrums of a symmetric $Z$-tensor and its absolute tensor}

In this section, we will study the relation between the spectrum of an even order symmetric $Z$-tensor
with nonnegative diagonal entries, and the spectrum of the absolute tensor of the $Z$-tensor.
It is proved that, if the symmetric $Z$-tensor is weakly irreducible and odd-bipartite, then the two spectral sets equal.
Furthermore, for an weakly irreducible symmetric $Z$-tensor with nonnegative diagonal enties,
we show that the spectral sets of the $Z$-tensor and its absolute tensor equal if and only if their spectral radii equal.
Before proving the conclusion, we firstly cite the definition of diagonal similar tensors, which is useful in the following analysis.

\bd\label{def61}$^{\cite{Shao13}}$ Let $\mathcal{A}$ and $\mathcal{B}$ be two order $m\geq2$ dimension $n$ tensors.
If there exists a nonsingular diagonal matrix $P$ of dimension $n$ such that $\mathcal{B}=P^{-(m-1)}\mathcal{A}P$,
then $\mathcal{A}$ and $\mathcal{B}$ are called diagonal similar.
\ed
Here, tensor $\mathcal{B}=P^{-(m-1)}\mathcal{A}P$ is defined by
$$b_{i_1i_2\cdots i_m}=\sum_{j_1,j_2,\cdots,j_m\in [n]}a_{j_1j_2\cdots j_m}p_{i_1j_1}^{m-1}p_{j_2i_2}\cdots p_{j_mi_m},~i_1,i_2,\cdots,i_m\in [n].$$

\bt\label{thm61} Assume order $m$ dimension $n$ symmetric $Z$-tensor $\mathcal{A}$ is defined as in (\ref{e51}). Suppose $\mathcal{C}$ is weakly irreducible.
Then, $\mathcal{A}$ and $|\mathcal{A}|$ are diagonal similar if and only if $m$ is even and $\mathcal{C}$ is weakly odd-bipartite.
\et
\proof For necessary, from Definition \ref{def61}, we know that there is a nonsingular diagonal matrix $P$ satisfying
$$\mathcal{A}=P^{-(m-1)}|\mathcal{A}|P,$$
i.e.,
$$\mathcal{D}-\mathcal{C}=P^{-(m-1)}(\mathcal{D}+\mathcal{C})P.$$
Since $\mathcal{D}=P^{-(m-1)}\mathcal{D}P$, we have
$$-\mathcal{C}=P^{-(m-1)}\mathcal{C}P,$$
which implies that
\begin{equation}\label{e61}
-c_{i_1i_2\cdots i_m}=c_{i_1i_2\cdots i_m}p_{i_1i_1}^{-(m-1)}p_{i_2i_2}\cdots p_{i_mi_m}.
\end{equation}
If $p_{11}=p_{22}=\cdots=p_{nn}$, by (\ref{e61}), we get $\mathcal{C}=0$, which is a
contradiction to the fact that $\mathcal{C}$ is weakly irreducible. So there are at least two distinct diagonal entries in $P$.

When $c_{i_1i_2\cdots i_m}\neq0$, by (\ref{e61}), one has
\begin{equation}\label{e62}
-p_{i_1i_1}^{m}=p_{i_1i_1}p_{i_2i_2}\cdots p_{i_mi_m}.
\end{equation}
By (\ref{e62}), and by the fact that $\mathcal{C}$ is weakly irreducible, we obtain
$$p_{ii}^m=p_{jj}^m,~~~i,j \in [n],$$
which implies that $m$ is even and
$$ V=\{i\in [n]~|~p_{ii}<0\}\neq \emptyset,\quad \tilde{V}=\{i\in [n]~|~p_{ii}>0\}\neq \emptyset.$$
Combining this with (\ref{e61})-(\ref{e62}), we know that
$$c_{i_1i_2\cdots i_m}=0,~~when ~|\{i_1,i_2,\cdots,i_m \}\cap V| ~ is ~even.$$
Thus, tensor $\mathcal{C}$ is weakly odd-bipartite corresponding to $V$ and the only if part holds.

For the if part, without loss of generality, suppose $\mathcal{C}$ is weakly odd-bipartite corresponding to $\Omega \subset [n]$.
Let $P$ be a diagonal matrix with $i$-th diagonal entries being -1 when $i\in \Omega$ and 1 when $i\notin \Omega$.
By a direct computation, one has
$$\mathcal{A}=P^{-(m-1)}|\mathcal{A}|P.$$
Apparently, $P$ is a nonsingular diagonal matrix. From Definition \ref{def61}, it follows that $\mathcal{A}$ and $|\mathcal{A}|$
are diagonal similar.\qed

It should be noted that diagonal similar tensors have
the same characteristic polynomials, and thus they have the same spectrum (see Theorem 2.1 of \cite{Shao13}), which is similar to the matrix case.

\bc\label{corol61}
Assume tensor $\mathcal{A}$ is defined as in Theorem \ref{thm61}. Let $m$ be even. Suppose $\mathcal{C}$ is odd-bipartite.
Then $Spec(\mathcal{A})=Spec(|\mathcal{A}|)$.
\ec

\bl\label{lema61}$^{\cite{YY11}}$ Let $\mathcal{A}$ and $\mathcal{B}$ be two order $m$ dimension $n$ tensors with $|\mathcal{B}|\leq \mathcal{A}$.
Then

(1) $\rho(\mathcal{B})\leq \rho(\mathcal{A})$.

(2) Furthermore, if A is weakly irreducible and $\rho(\mathcal{B})=\rho(\mathcal{A})$, where $\lambda=\rho(\mathcal{A})e^{i\psi}$ is an
eigenvalue of $\mathcal{B}$ with an eigenvector $y$, then,

(i) all the components of $y$ are nonzero;

(ii) let $U=diag(y_1/|y_1|,\cdots, y_n/|y_n|)$ be a nonsingular diagonal matrix, we have $\mathcal{B}=e^{i\psi}U^{-(m-1)}\mathcal{A}U$.
\el

\bt\label{thm62}  Assume order $m$ dimension $n$ symmetric $Z$-tensor $\mathcal{A}$ is defined as in (\ref{e51}). If
$\mathcal{C}$ is weakly irreducible, then $\rho(\mathcal{A})=\rho(|\mathcal{A}|)$ if and only if $Spec(\mathcal{A})=Spec(|\mathcal{A}|)$.
\et
\proof The sufficient condition is obvious. Now, we prove the only if part.
Suppose $\lambda=\rho(\mathcal{|\mathcal{A}|})e^{i\psi}$ is an eigenvalue of $\mathcal{A}$. Since $\mathcal{C}$ is weakly irreducible, from Lemma \ref{lema61}, we know that there exists a nonsingular diagonal matrix $P$ such that
\begin{equation}\label{e63}\mathcal{A}=e^{i\psi}P^{-(m-1)}|\mathcal{A}|P,
\end{equation}
which means
\begin{equation}\label{e64}
\mathcal{D}-\mathcal{C}=e^{i\psi}P^{-(m-1)}(\mathcal{D}+\mathcal{C})P.
\end{equation}
By the fact that all diagonal elements of $\mathcal{C}$ equal zero,
by (\ref{e64}), one has
$$\mathcal{D}=e^{i\psi}P^{-(m-1)}\mathcal{D}P=e^{i\psi}\mathcal{D},$$
which implies $e^{i\psi}=1$. So, by Definition \ref{def61} and (\ref{e63}), we know that $\mathcal{A}$ and $|\mathcal{A}|$ are
diagonal similar tensors. Thus, from Theorem 2.3 of \cite{Shao13}, it holds that $Spec(\mathcal{A})=Spect(|\mathcal{A}|)$.
\qed

\section{Final Remarks}
Odd-bipartite and even-bipartite tensors are defined in this paper. Using this,
we studied the relation between the largest H-eigenvalue of a $Z$-tensor with nonnegative diagonal elements, and
the largest H-eigenvalue of the $Z$-tensor's absolute tensor. Sufficient
and necessary conditions for the equality of these largest H-eigenvalues are given when the $Z$-tensor has even order. For the odd order case,
sufficient conditions are presented. Examples are given to verify the authenticity of the conclusions. On the other side, relation between spectral sets of an even order symmetric $Z$-tensor with nonnegative diagonal entries and its absolute tensor are studied.

In this paper, we only study the case of H-eigenvalues of $Z$-tensors. Do Z-eigenvalues of $Z$-tensors also hold in such case?
This may be an interesting work in the future.


\begin{thebibliography}{99}
\bibitem{CS13} D. Cartwright, B. Sturmfels, {\it The number of eigenvalues of a tensor}, Linear Algebra Appl. 438 (2013) 942-952.
\bibitem{Chang08} K.C. Chang, K. Pearson, T. Zhang,  {\it Perron Frobenius Theorem for nonnegative tensors}, Commu. Math. Sci. 6 (2008) 507-520.
\bibitem{CPZ09} K.C. Chang, K. Pearson, T. Zhang, {\it On eigenvalue problems of real symmetric tensors}, Journal of Mathematical Analysis and Applications 350 (2009) 416-422.
\bibitem{chen14} H. Chen, L. Qi, {\it Positive Definiteness and Semi-Definiteness of Even Order Symmetric Cauchy Tensors},
Journal of Industrial and Management Optimization 11 (2015) 1263-1274.
\bibitem{Chen14} Z. Chen, L. Qi, {\it Circulant Tensors with Applications to Spectral
Hypergraph Theory and Stochastic Process},  arXiv:1312.2752, 2014.
\bibitem{Cooper12} J. Cooper, A. Dutle, {\it Spectra of uniform hypergraphs}, Linear Algebra. Appl. 436 (2012) 3268-3292.
\bibitem{Ding13} W. Ding, L. Qi, Y. Wei, {\it M-Tensors and Nonsingular M-Tensors}, Linear Algebra Appl. 439 (2013) 3264-3278.
\bibitem{Fried13} S. Friedland, S. Gaubert, L. Han, {\it Perron-Frobenius theorem for nonnegative multilinear forms and extensions},
Linear Algebra Appl. 438 (2013) 738-749.
\bibitem{Qi2013} S. Hu, L. Qi, {\it The eigenvectors associated with the zero eigenvalues of the Laplacian and signless Laplacian tensors of a uniform}, Discrete Applied Mathematics  169 (2014) 140-151.
\bibitem{Hu13} S. Hu, L. Qi, J. Xie, {\it The largest Laplacian and signless Laplacian H-eigenvalues
of a uniform hypergraph}, Linear Algebra Appl. 469 (2015) 1-27.
\bibitem{Li13} G. Li, L. Qi, G. Yu, {\it Semismoothness of the maximum H-eigenvalue function of a symmetric tensor and its
application}, Linear Algebra Appl. 438 (2013) 212-218.
\bibitem{Lim05} L.-H. Lim, {\it Singular values and eigenvalues of tensors: a variational approach, in Proceedings
of the IEEE International Workshop on Computational Advances in Multi-
Sensor Addaptive Processing (CAMSAP¡¯05)}, 1 (2005) 129-132.
\bibitem{Mantica12} C.A. Mantica, L.G. Molinari, {\it Weakly Z-symmetric manifolds}, Acta Mathematica Hungarica 135(1) (2012) 80-96,
\bibitem{Ng09} M. Ng, L. Qi, G.L. Zhou, {\it Finding the largest eigenvalue of a nonnegative
tensor}, SIAM J. Matrx. Anal. Appl. 31 (2009) 1090-1099.
\bibitem{Oeding13} L. Oeding, G. Ottaviani, {\it Eigenvectors of tensors and algorithms for Waring decomposition}, Journal of Symbolic Computation 54 (2013) 9-35.
\bibitem{Pearson13} K. Pearson, T. Zhang, {\it On spectral hypergraph theory of the adjacency
tensor}, Graphs and Combinatorics 30 (2014) 1233-1248.
\bibitem{Qi05} L. Qi, {\it Eigenvalue of a real supersymmetric tensor}, J. Symb. Comput. 40 (2005) 1302-1324.
\bibitem{QiL13} L. Qi, {\it Symmetric nonnegative tensors and copositive tensors}, Linear Algebra Appl. 439 (2013) 228-238.
\bibitem{Qi13} L. Qi, {\it H$^+$-eigenvalues of Laplacian and signless Laplacian tensors}, Communications in Mathematical Sciences 12 (2014) 1045-1064.
\bibitem{Qi09} L. Qi, H.H. Dai, D. Han, {\it Conditions for strong ellipticity and M-eigenvalues}, Frontiers of
Mathmatics in China 4 (2009) 349-364.
\bibitem{Qi08}  L. Qi, Y. Wang, E. X. Wu, {\it D-eigenvalues of diffusion kurtosis tensor}, J. Comput. Appl. Math. 221 (2008) 150-157.
\bibitem{Shao13} J.Y. Shao, {\it A general product of tensors with applications}, Linear Algebra Appl. 439 (2013) 2350-2366.
\bibitem{Shao14} J.Y. Shao, H.Y. Shan, B. Wu, {\it Some Spectral Properties and Characterizations of Connected Odd-bipartite Uniform Hypergraphs}, to appear in: Linear and Multilinear Algebra. arXiv:1403.4845
\bibitem{Yang10} Y. Yang, Q. Yang, {\it Further results for Perron-Frobenius theorem for nonnegative tensors}, SIAM J.
Matrix Anal. Appl. 31 (2010) 2517-2530.
\bibitem{YY11} Y. Yang, Q. Yang, {\it On some properties of nonnegative weakly irreducible tensors}, arXiv preprint arXiv:1111.0713, 2011.
\bibitem{Zhang12} L. Zhang, L. Qi, G. Zhou, {\it M-tensors and some applications},  SIAM J. Matrix Anal. Appl. 35 (2014) 437-452.
\end{thebibliography}
\end{document}